\documentclass[12pt]{amsart}
\def\Int{\operatorname{Int}}
\def\diam{\operatorname{diam}}
\def\N{\mathbb{N}}

\newtheorem{Theorem}{Theorem}[section]

\newtheorem{Lemma}[Theorem]{Lemma}

\newtheorem{Thm}{Theorem}
\newtheorem{Cor}[Thm]{Corollary}
\theoremstyle{definition}

\theoremstyle{remark}

\begin{document}
\sloppy
\title{Parabolic subgroups of Coxeter groups acting by reflections on CAT(0) spaces}
\author{Tetsuya Hosaka} 
\address{Department of Mathematics, Faculty of Education, 
Utsunomiya University, Utsunomiya, 321-8505, Japan}
\email{hosaka@cc.utsunomiya-u.ac.jp}
\keywords{reflection group; Coxeter group; parabolic subgroup of a Coxeter group}
\subjclass[2000]{20F65; 20F55; 57M07}
\maketitle
\begin{abstract}
We consider a cocompact discrete reflection group $W$ of a CAT(0) space $X$.
Then $W$ becomes a Coxeter group.
In this paper, 
we study an analogy between the Davis-Moussong complex $\Sigma(W,S)$ and the CAT(0) space $X$, and 
show several analogous results 
about the limit set of a parabolic subgroup of the Coxeter group $W$.
\end{abstract}

\section{Introduction and preliminaries}

The purpose of this paper is to study 
the limit set of a parabolic subgroup of 
a reflection group of a CAT(0) space.
A metric space $(X,d)$ is called a {\it geodesic space} if 
for each $x,y \in X$, 
there exists an isometric embedding $\xi:[0,d(x,y)] \rightarrow X$ such that 
$\xi(0)=x$ and $\xi(d(x,y))=y$ (such a $\xi$ is called a {\it geodesic}).
We say that an isometry $r$ of a geodesic space $X$ is a {\it reflection} of $X$, 
if 
\begin{enumerate}
\item[(1)] $r^2$ is the identity of $X$, 
\item[(2)] $\Int F_r=\emptyset$ for the fixed-point set $F_r$ of $r$, 
\item[(3)] $X\setminus F_r$ has exactly two convex components $X_r^+$ and $X_r^-$, and 
\item[(4)] $r X_r^+=X_r^-$ and $r X_r^-=X_r^+$, 
\end{enumerate}
where the fixed-point set $F_r$ of $r$ is called the {\it wall} of $r$.
Let $X_r^+$ and $X_r^-$ be the two convex connected components of $X\setminus F_r$, 
where $X_r^+$ contains a basepoint of $X$.
An isometry group $\Gamma$ of a geodesic space $X$ 
is called a {\it reflection group}, 
if some set of reflections of $X$ generates $\Gamma$.

Let $\Gamma$ be a reflection group of a geodesic space $X$ and 
let $R$ be the set of all reflections of $X$ in $\Gamma$.
Now we suppose that 
the action of $\Gamma$ on $X$ is proper, 
that is, $\{\gamma\in\Gamma\,|\, \gamma x\in B(x,N)\}$ is finite 
for any $x\in X$ and $N>0$ (cf.\ \cite[p.131]{BH}).
Then the set $\{F_r\,|\, r\in R\}$ is locally finite.
Let $C$ be a component of $X\setminus \bigcup_{r\in R} F_r$, 
which is called a {\it chamber}.
Then 
$\Gamma C=X\setminus \bigcup_{r\in R} F_r$, 
$\Gamma \overline{C}=X$ 
and for each $\gamma\in\Gamma$, 
either $C\cap \gamma C=\emptyset$ or $C=\gamma C$.
We say that $\Gamma$ is a {\it cocompact discrete reflection group} of $X$, 
if $\overline{C}$ is compact and $\{\gamma\in\Gamma\,|\,C=\gamma C\}=\{1\}$.
Every Coxeter group is a cocompact discrete reflection group 
of some CAT(0) space.

A {\it Coxeter group} is a group $W$ having a presentation
$$\langle \,S \, | \, (st)^{m(s,t)}=1 \ \text{for}\ s,t \in S \,
\rangle,$$ 
where $S$ is a finite set and 
$m:S \times S \rightarrow \N \cup \{\infty\}$ is a function 
satisfying the following conditions:
\begin{enumerate}
\item[(1)] $m(s,t)=m(t,s)$ for any $s,t \in S$,
\item[(2)] $m(s,s)=1$ for any $s \in S$, and
\item[(3)] $m(s,t) \ge 2$ for any $s,t \in S$ such that $s\neq t$.
\end{enumerate}
The pair $(W,S)$ is called a {\it Coxeter system}.
H.S.M.~Coxeter showed that 
a group $\Gamma$ is a finite reflection group of some Euclidean space 
if and only if $\Gamma$ is a finite Coxeter group. 
Every Coxeter system $(W,S)$ induces 
the Davis-Moussong complex $\Sigma(W,S)$ which is a CAT(0) space 
(\cite{D1}, \cite{D2}, \cite{M}).
Then 
the Coxeter group $W$ is a cocompact discrete reflection group 
of the CAT(0) space $\Sigma(W,S)$.
It is known that 
a group $\Gamma$ is a cocompact discrete reflection group 
of some geodesic space if and only if 
$\Gamma$ is a Coxeter group (\cite{H5}).

Let $W$ be a cocompact discrete reflection group of a CAT(0) space $X$, 
let $R$ be the set of reflections in $W$, 
let $C$ be a chamber and 
let $S$ be a {\it minimal} subset of $R$ such that $C=\bigcap_{s\in S}X_s^+$ 
(i.e.\ $C\neq \bigcap_{s\in S\setminus\{s_0\}}X_s^+$ for any $s_0\in S$).
Then $\langle S\rangle C=X\setminus \bigcup_{r\in R} F_r=WC$, 
$S$ generates $W$ and 
the pair $(W,S)$ is a Coxeter system (\cite{H5}).
For a subset $T$ of $S$, 
$W_T$ is defined as the subgroup of $W$ generated by $T$, 
and called a {\it parabolic subgroup}.
It is known that 
the pair $(W_T,T)$ is also a Coxeter system.

Let $X$ be a CAT(0) space and let $\Gamma$ be a group 
which acts properly by isometries on $X$.

The {\it limit set of $\Gamma$} ({\it with respect to $X$})
is defined as 
$$ L(\Gamma) = \overline{\Gamma x_0} \cap\partial X, $$
where $\overline{\Gamma x_0}$ is the closure of the orbit $\Gamma x_0$ in $X \cup \partial X$ 
and $x_0$ is a point in $X$.
We note that the limit set $L(\Gamma)$ 
is independent of the point $x_0 \in X$.

Also we say that (the action of) $\Gamma$ is 
{\it convex-cocompact}, 
if there exists a compact subset $K$ of $X$ such that 
${\mathcal R}_{x_0}(L(\Gamma)) \subset \Gamma K$ for some $x_0\in X$, 
where ${\mathcal R}_{x_0}(L(\Gamma))$ is the union of the images of 
all geodesic rays $\xi$ issuing from $x_0$ with $\xi(\infty) \in L(\Gamma)$.
We note that 
for a group acting on a proper CAT(0) space, 
``convex-cocompactness'' agrees with 
``geometrically finiteness'' (cf.\ \cite{H2} and \cite{H3}).

We first prove the following theorem in Section~2.

\begin{Thm}
For each subset $T\subset S$, 
\begin{enumerate}
\item[(1)] $W_T \overline{C}$ is convex (hence CAT(0)),
\item[(2)] the limit set $L(W_T)$ of $W_T$ 
coincides with the boundary $\partial (W_T \overline{C})$, and 
\item[(3)] the action of $W_T$ on $X$ is convex-cocompact.
\end{enumerate}
\end{Thm}

This theorem implies the following corollary (\cite{H2} and \cite{H3}).

\begin{Cor}
For each subset $T\subset S$, 
the following statements are equivalent:
\begin{enumerate}
\item[(1)] $[W:W_T]<\infty$;
\item[(2)] $L(W_T)=\partial X$;
\item[(3)] $\Int_{\partial X}L(W_T)\neq\emptyset$.
\end{enumerate}
\end{Cor}

In Section~3, 
we show the following theorem which is an analogue of Lemma~4.2 in \cite{H1}.

\begin{Thm}
Let $x_0\in C$ and let $w\in W$.
Then there exists a reduced representation 
$w=s_1\cdots s_l$ such that 
$$d_H([x_0,wx_0],P_{s_1,\dots,s_l})\le \diam \overline{C},$$ 
where $d_H$ is the Hausdorff distance and 
$P_{s_1,\ldots,s_l}=
[x_0,s_1x_0] \cup [s_1x_0,(s_1s_2)x_0] \cup \cdots \cup
[(s_1\cdots s_{l-1})x_0,wx_0]$.
\end{Thm}

Using this theorem,
we can obtain the following corollaries 
by the same argument used in \cite{H1} and \cite{H4}.

\begin{Cor}
For each subset $T\subset S$, 
the limit set $L(W_T)$ is $W$-invariant if and only if 
$W=W_{\tilde{T}}\times W_{S\setminus \tilde{T}}$.
\end{Cor}

Here $W_{\tilde{T}}$ is the {\it essential parabolic subgroup of} $W_T$
(cf.\ \cite{H1}), that is, 
$W_{\tilde{T}}$ is the minimum parabolic subgroup 
of finite index in $(W_T,T)$.

We denote by $o(g)$ the order of 
an element $g$ in the Coxeter group $W$.
For $s_0\in S$, we define 
$W^{\{s_0\}}=\{w\in W\,|\,\ell(ws)>\ell(w) 
\ \text{for each}\ s\in S\setminus\{s_0\} \}\setminus \{1\}$.
A subset $A$ of a space $Y$ is said to be {\it dense} in $Y$, 
if $\overline{A}=Y$. 
A subset $A$ of a metric space $Y$ is said to be {\it quasi-dense}, 
if there exists $N>0$ such that 
each point of $Y$ is $N$-close to some point of $A$.

\begin{Cor}
Suppose that 
$W^{\{s_0\}}$ is quasi-dense in $W$ 
with respect to the word metric and 
$o(s_0t_0)=\infty$ for some $s_0,t_0\in S$.
Then there exists $\alpha \in \partial X$ 
such that the orbit $W\alpha$ is dense in $\partial X$.
\end{Cor}

\begin{Cor}
If the set 
$$\bigcup\{W^{\{s\}}\,|\, s\in S \ \text{such that}\ o(st)=\infty \ \text{for some}\ t\in S\}$$ 
is quasi-dense in $W$, 
then $\{w^\infty\,|\, w\in W \ \text{such that}\ o(w)=\infty\}$ 
is dense in $\partial X$.
\end{Cor}

A subset $T$ of $S$ is said to be {\it spherical}, if $W_T$ is finite.

\begin{Cor}
Suppose that there exist a maximal spherical subset $T$ of $S$ 
and an element $s_0\in S$ such that $o(s_0t)\ge 3$ for each $t\in T$
and $o(s_0t_0)=\infty$ for some $t_0\in T$. 
Then 
\begin{enumerate}
\item[(1)] $W\alpha$ is dense in $\partial X$ 
for some $\alpha \in \partial X$, and 
\item[(2)] 
$\{w^\infty\,|\, w\in W \ \text{such that}\ o(w)=\infty\}$ 
is dense in $\partial X$.
\end{enumerate}
\end{Cor}

\section{Convex-cocompactness of parabolic subgroups}

Let $W$ be a cocompact discrete reflection group of a CAT(0) space $X$, 
let $C$ be a chamber containing a basepoint of $X$,
let $R$ be the set of reflections in $W$, and 
let $S$ be a minimal subset of $R$ such that $C=\bigcap_{s\in S}X_s^+$.
(Then the pair $(W,S)$ is a Coxeter system \cite{H5}.) 
For each reflection $r$ in $W$, 
$F_r$ is the wall of $r$ and 
$X_r^+$ and $X_r^-$ are the two convex components of $X\setminus F_r$ 
such that $C\subset X_r^+$ and $C\cap X_r^-=\emptyset$.
We note that $F_r$, $X_r^+\cup F_r$ and $X_r^-\cup F_r$ are convex.

The following lemmas are known.

\begin{Lemma}[{\cite[Lemma~3.4]{H5}}]\label{lem0}
Let $w\in W$ and let $s\in S$.
Then $\ell(w)<\ell(sw)$ if and only if 
$wC\subset X^+_s$.
\end{Lemma}

\begin{Lemma}[{\cite[Lemma~1.3]{D3}}]\label{lem0-1}
Let $w\in W$ and let $T\subset S$.
Then there exists a unique element of 
shortest length in the coset $W_Tw$.
Moreover, the following statements are equivalent:
\begin{enumerate}
\item[(1)] $w$ is the element of shortest length in the coset $W_Tw$;
\item[(2)] $\ell(sw)>\ell(w)$ for any $s\in T$;
\item[(3)] $\ell(vw)=\ell(v)+\ell(w)$ for any $v\in W_T$.
\end{enumerate}
\end{Lemma}

We first show the following lemma.

\begin{Lemma}\label{lem1}
Let $T\subset S$.
Then $w X^+_s=X^+_{wsw^{-1}}$
for any $w\in W_T$ and $s\in S\setminus T$.
\end{Lemma}

\begin{proof}
Let $T\subset S$, $w\in W_T$ and $s\in S\setminus T$.
Then $\ell(sw^{-1})>\ell(w^{-1})$.
Hence $w^{-1}C\subset X^+_s$ by Lemma~\ref{lem0}.
Thus $C\subset w X^+_s$, i.e., 
$w X^+_s=X^+_{wsw^{-1}}$.
\end{proof}

Using lemmas above, we prove the following theorem.

\begin{Theorem}\label{thm1}
For each subset $T\subset S$, 
\begin{enumerate}
\item[(1)] $W_T \overline{C}$ is convex (hence CAT(0)),
\item[(2)] the limit set $L(W_T)$ of $W_T$ 
coincides with the boundary $\partial (W_T \overline{C})$, and 
\item[(3)] the action of $W_T$ on $X$ is convex-cocompact.
\end{enumerate}
\end{Theorem}

\begin{proof}
Let $T\subset S$. 
Then we show that 
$$W_T \overline{C} 
=\bigcap\{\overline{X^+_{wsw^{-1}}}\,|\, w\in W_T,\;s\in S\setminus T\}.$$

For each $v,w \in W_T$ and $s\in S\setminus T$, 
$C\subset v^{-1}w X^+_s$ by Lemma~\ref{lem1}.
Hence $vC\subset wX^+_s=X^+_{wsw^{-1}}$ by Lemma~\ref{lem1}.
Thus $v\overline{C}\subset \overline X^+_{wsw^{-1}}$ 
for any $v,w \in W_T$ and $s\in S\setminus T$, 
that is, 
$$W_T \overline{C} \subset
\bigcap\{\overline{X^+_{wsw^{-1}}}\,|\, w\in W_T,\;s\in S\setminus T\}.$$

To prove 
$$W_T \overline{C} \supset
\bigcap\{\overline{X^+_{wsw^{-1}}}\,|\, w\in W_T,\;s\in S\setminus T\},$$
we show that for each $v\in W\setminus W_T$, 
there exist $w\in W_T$ and $s\in S\setminus T$ 
such that $vC\subset X^-_{wsw^{-1}}$.
Let $v\in W\setminus W_T$.
By Lemma~\ref{lem0-1}, 
There exists a unique element $x\in W_T v$ of shortest length.
Let $w=vx^{-1}$.
Here we note that $w\in W_T$ and $\ell(v)=\ell(w)+\ell(x)$.
Let $s\in S$ such that $\ell(sx)<\ell(x)$.
By Lemma~\ref{lem0-1}~(2), $s\in S\setminus T$.
Then 
$$ \ell(sw^{-1}v)=\ell(sx)<\ell(x)=\ell(w^{-1}v). $$
Hence $w^{-1}vC\subset X^-_s$ by Lemma~\ref{lem0}.
By Lemma~\ref{lem1}, 
$vC\subset wX^-_s=X^-_{wsw^{-1}}$.
Therefore
$$W_T \overline{C} 
=\bigcap\{\overline{X^+_{wsw^{-1}}}\,|\, w\in W_T,\;s\in S\setminus T\}.$$

Since $\overline{X^+_{wsw^{-1}}}=X^+_{wsw^{-1}}\cup F_{wsw^{-1}}$ is convex 
for any $w\in W_T$ and $s\in S\setminus T$, 
$W_T \overline{C}$ is convex.
Hence 
$L(W_T)=\partial (W_T \overline{C})$ and 
the action of $W_T$ on $X$ is convex-cocompact.
\end{proof}

\section{On geodesics and reduced representations}

We give the following lemma 
which is an analogue of a result about Davis-Moussong complexes.

\begin{Lemma}\label{lem3-1}
Let $w\in W$, let $w=s_1\cdots s_l$ be a reduced representation 
and let $T=\{s_1,\dots,s_l\}$.
Then 
$$
\overline{C}\cap w\overline{C}
=\bigcap_{t\in T}(F_t\cap \overline{C}) 
=\bigcap_{t\in T}(t\overline{C}\cap \overline{C}) 
=\bigcap_{v\in W_T}v\overline{C}.
$$
\end{Lemma}

\begin{proof}
Let $y\in \overline{C}\cap w\overline{C}$.
Since $\ell(s_1 w)<\ell(w)$, 
$wC\subset X^-_{s_1}$ by Lemma~\ref{lem0}.
Then
$$y\in \overline{C}\cap w\overline{C}\subset 
\overline{X^+_{s_1}}\cap \overline{X^-_{s_1}}=F_{s_1}.$$
Hence $s_1y=y$ and 
$$y=s_1y\in s_1(\overline{C}\cap w\overline{C})
=s_1\overline{C}\cap (s_2\cdots s_l)\overline{C},$$
i.e., $y\in\overline{C}\cap (s_2\cdots s_l)\overline{C}$.
By iterating the above argument, 
$s_iy=y$ for any $i\in\{1,\dots,l\}$, that is, 
$ty=y$ for any $t\in T$.
Hence $y\in \bigcap_{t\in T}(F_t\cap \overline{C})$.
Thus 
$\overline{C}\cap w\overline{C}
\subset\bigcap_{t\in T}(F_t\cap \overline{C})$.

Since $F_t\cap \overline{C}=t\overline{C}\cap\overline{C}$ 
for any $t\in T$, 
$\bigcap_{t\in T}(F_t\cap \overline{C})
=\bigcap_{t\in T}(t\overline{C}\cap \overline{C})$.

Let $y\in\bigcap_{t\in T}(F_t\cap \overline{C})$.
Then $ty=y$ for any $t\in T$.
Since $T$ generates $W_T$, $vy=y$ for any $v\in W_T$.
Hence $y=vy\in v\overline{C}$ for each $v\in W_T$.
Thus 
$\bigcap_{t\in T}(F_t\cap \overline{C})
\subset\bigcap_{v\in W_T}v\overline{C}$.

It is obvious that 
$\bigcap_{v\in W_T}v\overline{C}\subset
\overline{C}\cap w\overline{C}$, 
since $1,w\in W_T$.
\end{proof}

\begin{Lemma}\label{lem3-2}
Let $w\in W$, let $w=s_1\cdots s_l$ be a reduced representation 
and let $T=\{s_1,\dots,s_l\}$.
Then 
$\overline{C}\cap w\overline{C}\neq\emptyset$ 
if and only if 
$W_T$ is finite.
\end{Lemma}

\begin{proof}
Suppose that $\overline{C}\cap w\overline{C}\neq\emptyset$.
Then 
$\bigcap_{v\in W_T}v\overline{C}\neq\emptyset$ by Lemma~\ref{lem3-1}.
Hence $W_T$ is finite 
because the action of $W$ on $X$ is proper.

Suppose that $W_T$ is finite.
Then 
$W_T$ acts on the CAT(0) space $W_T\overline{C}$ 
by Theorem~\ref{thm1}.
By \cite[Corollary~II.2.8(1)]{BH}, 
there exists a fixed-point $y\in W_T\overline{C}$ 
such that $vy=y$ for any $v\in W_T$.
Then $y\in \overline{C}\cap w\overline{C}$ 
which is non-empty.
\end{proof}

By the proof of \cite[Lemma~4.2]{H1}, 
we can obtain the following theorem 
from Lemmas~\ref{lem0}, \ref{lem3-1} and \ref{lem3-2}.

\begin{Theorem}
Let $x_0\in C$ and let $w\in W$.
Then there exists a reduced representation
$w=s_1\cdots s_l$ such that 
$$d_H([x_0,wx_0],P_{s_1.\dots,s_l})\le \diam \overline{C},$$ 
where $d_H$ is the Hausdorff distance and 
$P_{s_1,\ldots,s_l}=
[x_0,s_1x_0] \cup [s_1x_0,(s_1s_2)x_0] \cup \cdots \cup
[(s_1\cdots s_{l-1})x_0,wx_0]$.
\end{Theorem}

%

%
\end{document}